\documentclass[11pt,a4paper,leqno]{article}
\usepackage{amsmath}
\usepackage{amssymb}
\advance\topmargin by -2.2cm
\newcommand{\bbd}{I\!\!D}
\newcommand{\bbr}{I\!\!R}
\newcommand{\bbn}{I\!\!N}
\newcommand{\bbz}{Z\!\!\!Z}
\newcommand{\bbf}{I\!\!F}

\newcommand{\calb}{{\cal B}}
\newcommand{\calc}{{\cal C}}

\newcommand{\cale}{{\cal E}}
\newcommand{\calf}{{\cal F}}

\newcommand{\cali}{{\cal I}}

\newcommand{\calh}{{\cal H}}

\newcommand{\call}{{\cal L}}

\newcommand{\caln}{{\cal N}}

\newcommand{\barr}{\begin{array}}
\newcommand{\earr}{\end{array}}
\newcommand{\beqq}{\begin{equation}}
\newcommand{\eeqq}{\end{equation}}
\newcommand{\beao}{\begin{eqnarray*}}
\newcommand{\eeao}{\end{eqnarray*}\noindent}
\newcommand{\beam}{\begin{eqnarray}}
\newcommand{\eeam}{\end{eqnarray}\noindent}

\newcommand{\halmos}{\quad\hfill\mbox{$\Box$}}

\newcommand{\la}{\lambda}

\newcommand{\si}{\sigma}
\newcommand{\al}{\alpha}
\newcommand{\vth}{\vartheta}

\newcommand{\wh}{\widehat}
\newcommand{\wt}{\widetilde}

\newcommand{\lra}{\longrightarrow}

\newcommand{\nto}{n\to\infty}

\begin{document}

{\huge 
Estimating a periodicity parameter\\ in the drift of a time inhomogeneous diffusion}
\vskip0.5cm
{\bf R.\ H\"opfner}, Johannes Gutenberg Universit\"at Mainz\\
{\bf Yu.\ Kutoyants}, Universit\'e du Maine, Le Mans\\


{\bf Abstract: } We consider a diffusion $(\xi_t)_{t\ge 0}$ whose drift contains some deterministic periodic signal. Its shape being fixed and known, up to scaling in time, the periodicity of the signal is the unknown parameter $\vth$ of interest. 
We consider sequences of local models at $\vth$, corresponding to continuous observation of the process $\xi$ on the time interval $[0,n]$ as $n\to\infty$, with suitable choice of local scale at $\vth$. Our tools --under an ergodicity condition-- are path segments of $\xi$  corresponding to the period $\vth$, and limit theorems for certain functionals of the process $\xi$ which are not additive functionals. 
\\
When the signal is smooth, with local scale $n^{-3/2}$ at $\vth$, we have local asymptotic normality (LAN) in the sense of Le Cam (1969).  
When the signal has a finite number of discontinuities, with local scale $n^{-2}$ at $\vth$, we obtain a limit experiment of different type, studied by Ibragimov and Khasminskii (1981), where smoothness of the parametrization (in the sense of Hellinger distance) is H\"older $\frac12$.
{\bf MSC: \quad 62 F 12 , 60 J 60} \\


Write $P^\vth$ for the law on $\,\left( C([0,\infty)),\calc([0,\infty)) \right)\,$ of the solution to the SDE 
\beqq\label{process}
d\xi^\vth_t  \;=\;  \left[ S(\vth,t) + b(\xi^\vth_t) \right]dt   \;+\;  \si(\xi^\vth_t)\,dW_t  \quad,\quad t\ge 0 \;\;,\;\; \xi^\vth_0=x_0  
\eeqq
whose drift involves a periodic deterministic signal 
$t \to S(\vth,t)$
governed by an unknown parameter  
\beqq\label{parametrization}
S(\vth,t) \;:=\; S_0(\; \frac1\vth\, t\, )   \quad,\quad \vth\in\Theta:=(0,\infty) \;,\; t\ge 0
\eeqq
where the $1$-periodic function $S_0(\cdot)$ is fixed and known. By (\ref{parametrization}), the signal $t\to S(\vth,t)$  contained in the drift of the process (\ref{process}) is $\vth$-periodic under $P^\vth$.

A special case of our problem (\ref{process})+(\ref{parametrization}) is the  'signal in white noise' model $b(\cdot)\equiv 0$, $\si(\cdot)\equiv 1$ for which Ibragimov and Khasminskii [IH 81, pp.\ 209] proved local asymptotic normality at rate $n^{-3/2}$ when observing the process under unknown $\vth$ up to time $n$, provided the signal $S_0(\cdot)$ is smooth. Their approach was combined with $L^2$-methods by Golubev [G 88] to estimate the period $\vth$ together with the shape of the signal, at the same rate $n^{-3/2}$. Based on this, estimation of the shape of the signal under unknown periodicity was considered by Castillo, L\'evy-Leduc and Matias [CLM 06] who prove nonparametric rates over nonparametric function classes. We will not go in the last mentioned  directions. For discontinuous signals, the features are essentially different. Local scale is $n^{-2}$, and a different type of limit experiment arises, also studied in [IH 81]: likelihoods between laws $\wt P_h$ and $\wt P_0$ in the limit model $\{ \wt P_h : h\in\bbr \}$ are of type $\exp\{ \wt W(h) -\frac12|h|\}$ for double-sided Brownian motion $(\wt W(h))_{h\in\bbr}$. This limit model has a number of interesting properties: in particular, for quadratic loss, the Bayes estimator in the limit experiment is better than the maximum likelihood estimator ([IH 81, p. 342], [RS 95]). Note that for the 'signal in white noise' setting $b(\cdot)\equiv 0$, $\si(\cdot)\equiv 1$, the observed process is always a Gaussian process, hence specific techniques for Gaussian processes are available which do not carry over to the problem (\ref{process})+(\ref{parametrization}).

To illustrate the type of difficulty which arises with nontrivial drift and diffusion coefficient, we anticipate and mention one typical problem  of convergence which arises when we consider asymptotics of likelihoods. Take the example of discontinuous $S_0(\cdot)$ defined from $1$-periodic continuation of 
$$
S_0(t) \;:=\; 1_{]r_1,r_2]}(t) \quad\mbox{for $0\le t\le 1$}
$$
with given jump times $0<r_1<r_2<1$. Based on continuous observation of (\ref{process}) under unknown $\vth$ up to time $n$, the local structure of log-likelihoods close to $\vth$, with local scale $n^{-2}$ and local parameter $h$ at $\vth$, makes appear objects like 
$$
\sum_{k=0}^{\lfloor \frac{n}{\vth}\rfloor } \int_{(\vth+h n^{-2})(k+r_i)}^{\vth(k+r_i)} \frac{1}{\si^2(\xi^\vth_s)}\, ds  
\quad\mbox{if $h<0$} \quad,\quad 
\sum_{k=0}^{\lfloor \frac{n}{\vth}\rfloor } \int_{\vth(k+r_i)}^{(\vth+h n^{-2})(k+r_i)} \frac{1}{\si^2(\xi^\vth_s)}\, ds 
\quad\mbox{if $h>0$} \;. 
$$
If $\si(\cdot)\equiv 1$ is constant, there is nothing to prove, the above quantities being deterministic and close to $|h|\cdot\frac{1}{2\vth^2}$. For nonconstant $\si(\cdot)$, proving that indeed the above quantities converge in $P^\vth$-probability as $\nto$ to some deterministic limit which contains a factor $|h|$, under suitable assumptions on the process under $\vth$, requires quite some work: the proofs will have to exploit the periodicity structure of the process in order to work with 'ergodic properties', will have to control fluctations of the process $\xi^\vth$ over the above intervals of integration $[(\vth+h n^{-2})(k+r_i) , \vth(k+r_i)]$ or $[\vth(k+r_i) , (\vth+h n^{-2})(k+r_i)]$ whose length is $\approx\frac{1}{n^2}k$, and so on. Hence, proving convergence of likelihood ratios in the problem (\ref{process})+(\ref{parametrization}) is very different from what one does in the 'signal in white noise' model.

In this paper, we will consider sequences of local models at $\vth$ in the problem (\ref{process})+(\ref{parametrization}), corresponding to continuous observation of the process (\ref{process}) under $\vth$ on the time interval $[0,n]$ as $n\to\infty$, with suitable choice of local scale at $\vth$. Our tools --under an ergodicity condition-- are path segments of $\xi$  corresponding to the period $\vth$ which form a time homogeneous Markov chain of path segments, and limit theorems for certain functionals of the process $\xi$ which are not additive functionals.

When the signal is smooth, with local scale $n^{-3/2}$ at $\vth$, we will prove in theorem 1.1 below local asymptotic normality (LAN) in the sense of Le Cam [L 69] (studied intensely and under a broad range of aspects by many authors, see e.g.\ Hajek and Sidak [HS 67], Hajek [H 70], Ibragimov and Khasminskii [IH 81], Strasser [85], Davies [85], Le Cam and Yang [LY 90], Pfanzagl [P 94], Witting and M\"uller-Funk [WM 95], Kutoyants [K 98], van der Vaart [V 98], Kutoyants [K 04], Liese and Miescke~[08]). In the limit experiment, we have one estimator for the local parameter --the central statistics-- which minimizes the risk simultaneously under a broad class of loss functions.

When the signal is discontinuous, with local scale $n^{-2}$ at $\vth$, we will prove in theorem 1.2 below convergence of the log-likelihoods in local models at $\vth$ to a limit experiment of type 
$$
\left\{ \wt P_h : h\in\bbr \right\} \quad\mbox{with likelihoods}\quad \wt L^{h / 0} \;=\; \exp\left\{ \wt W(h) -\frac12|h| \right\}
$$
where $(\wt W(h))_{h\in\bbr}$ is double-sided Brownian motion. This limit experiment, studied by Ibragimov and Khasminskii for discontinuous signals in white noise [IH 81, sections VII.2--3] and in the sequel by several authors in various contexts (see e.g.\  Golubev [G 79], Deshayes and Picard [DP 84], Rubin and Song [RS 95], K\"uchler and Kutoyants [KK 00], Kutoyants [K 04, section 3.4], Dachian [D 10], H\"opfner and Kutoyants [HK 10]),  is very different from the previous one: the smoothness of the likelihood as a function of the parameter $h$ is the smoothness of the Brownian path, the smoothness of parametrization in the sense of Hellinger distance is H\"older $\frac12$ (i.e.:  Hellinger distances between $\wt P_{h'}$ and $\wt P_h$  are of order $|h'-h|^{1/2}$ as $h'\to h$), only loss-function-specific minimax results are known, and a squared loss Bayes estimator is strictly better than the maximum likelihood estimator (where better means: with respect to squared loss); no tool is known which permits to compare the minimax bounds related to  different loss functions.

\section*{1. Main results}

For the problem (\ref{process})+(\ref{parametrization}), we assume Lipschitz conditions on $b(\cdot)$ and $\si(\cdot)$, and have (e.g.~[KS 91]) for every $\vth\in\Theta$ a unique strong solution of (\ref{process}). The starting point $x_0$ in (\ref{process}) is fixed and does not vary with $\vth$. Our main assumption on the process will be 
$$
\mbox{ {\bf (H2)} :}\quad \mbox{for every $\vth\in\Theta$, the $\vth$-grid chain $(\xi^\vth_{k\vth})_k$ under $\vth$ is positive Harris recurrent} \;. 
$$
As an example, condition {\bf (H2)} holds for piecewise continuous $S_0(\cdot)$ in case where 
the process without signal is an ergodic Ornstein Uhlenbeck process ([HK 10, example 2.3], [DFK 10]). Under {\bf (H2)}, we write $\mu^\vth$ for the (unique) invariant law of the $\vth$-grid chain $(\xi^\vth_{k\vth})_k$ under $\vth$. We will work under either 'smooth signal' hypotheses 
\beao 
\mbox{ {\bf (H0)} :}\quad &&\mbox{the $1$-periodic function $S_0(\cdot)$ in (\ref{parametrization}) is $\calc^2$ on $[0,\infty)$}   \\
\mbox{ {\bf (H1)} :}\quad &&\mbox{$\si(\cdot)$ is bounded away from $0$}    
\eeao 
or 'discontinuous signal' hypotheses 
\beao 
\mbox{ {\bf (H0')} :}\quad &&\mbox{$S_0(\cdot)$ is $1$-periodic, piecewise continuous, and Lipschitz between the jumps}   \\
\mbox{ {\bf (H1')} :}\quad &&\mbox{$\si(\cdot)$ is bounded away from both $0$ and $\infty$}  
\eeao 
and stress that {\bf (H0')} deals with discontinuous functions $S_0(\cdot)$ with finitely many jumps on $(0,1]$. \\

Let $\,\left( C([0,\infty)),\calc([0,\infty)), \bbf \right)\,$ denote the canonical path space for solutions of (\ref{process}) and 
$\eta=(\eta_t)_{t\ge 0}$ the canonical process; i.e.\ $\bbf =(\calf_t)_{t\ge 0}$ where $\calf_t$ is generated by observation of $\eta$ up to time $t+$. $\;(P^\vth_{s,t})_{0\le s< t<\infty}$ denotes the semigroup of transition probabilities of the process (\ref{process}) under $\vth$. Write $P^\vth_t$ for the restriction of $P^\vth$ to $\calf_t$. Let $L^{\zeta / \vth}=(L^{\zeta / \vth}_t)_{t\ge 0}$ denote the likelihood ratio process of $P^\zeta$ to $P^\vth$ relative to the filtration $\bbf$ (cf.\ [LS 81], [JS 87], [K 04]): 
\beqq\label{likelihoodratioprocess}
L^{\zeta / \vth} 
\;=\;  \cale_\vth\left( \int_0^\cdot \frac{S(\zeta,t)-S(\vth,t)}{\si^2(\eta_t)}\; dm^{(\vth)}_t \right) 
\;=\;  \cale_\vth\left( \int_0^\cdot \frac{S(\zeta,t)-S(\vth,t)}{\si(\eta_t)}\; dB^{(\vth)}_t \right) 
\eeqq
where $m^{(\vth)}$ is the martingale part of the canonical process under $P^\vth$, and where $B^{(\vth)} := \int_0^\cdot \frac{1}{\si(\eta_s)}dm^{(\vth)}_s$ is a  $(\bbf, P^\vth)$-Brownian motion. We shall consider the sequence of experiments 
\beqq\label{sequence_of_experiments}
\left(\, C([0,\infty)) \,,\, \calf_{n} \,,\, \left\{ P^{\zeta}_{n} : \zeta\in\Theta \right\} \,\right) \quad,\quad n\to\infty 
\eeqq
locally in small neighbourhoods of some fixed  $\vth\in\Theta$.

Under 'smooth signal' hypotheses, the limit of local models (\ref{Gaussianlocalmodels}) at $\vth$, with local scale $n^{-3/2}$ (thus essentially faster than the usual $n^{-1/2}$ in time homogenous ergodic diffusions, and also essentially faster than the rate $n^{-1/2}$ in [HK 11] when the drift contains a parametrized continuous signal of known periodicity), will be the well known Gaussian shift model. The following is a 2nd Le Cam lemma (in the language of Hajek and Sidak [HS 67]) for estimation of the periodicity in the problem (\ref{process})+(\ref{parametrization}), for smooth signals.\\

{\bf 1.1 Theorem : } Under hypotheses {\bf (H0)}+{\bf (H1)}+{\bf (H2)} consider the sequence of experiments (\ref{sequence_of_experiments}) and local models at $\vth\in\Theta$
\beqq\label{Gaussianlocalmodels}
\left\{ P^{\vth + n^{-3/2} h}_{n} : h\in\Theta_{\vth,n} \right\}  \quad,\quad n\to\infty 
\eeqq
where  $\Theta_{\vth,n}:= \{ h\in\bbr: \vth + n^{-3/2} h \in\Theta\}$. 

a) For every $\vth\in \Theta$, we have LAN at $\vth$ with local scale $n^{-3/2}$ and Fisher information 
\beqq\label{fisherinformation} 
\cali^{(\vth)} \;\;=\;\;  
\frac{1}{ 3\, \vth^4 }\, \int_0^1 dv\; [S_0']^2(v)\; [\mu^\vth P^\vth_{0,v\vth}]( \frac{1}{\si^2} ) \;. 
\eeqq

b) For every $\vth\in\Theta$, for arbitrary bounded sequences $(h_n)_n$ in $\bbr^d$, we have a quadratic decomposition of log-likelihood ratios  
\beqq\label{quadraticdecomposition}
\log L_{n}^{(\vth + n^{-3/2} h_n) / \vth}  \;=\; h_n\, \Delta_n^{(\vth)} \;-\; \frac12 h_n^2\, \cali^{(\vth)} \;+\; o_{P^\vth}(1) \quad\mbox{as $\nto$} 
\eeqq
with score 
\beqq\label{score}
\Delta_n^{(\vth)} \;=\; \frac{1}{\sqrt{n^3}} \int_0^{n} \frac{\dot S(\vth,s)}{\si(\eta_s)}\, dB^{(\vth)}_s 
\;=\;  \frac{-\vth^{-2}}{\sqrt{n^3}} \int_0^{n} \frac{s\; S_0'(\,\frac{1}{\vth}\, s\,)}{\si(\eta_s)}\, dB^{(\vth)}_s 
\eeqq
such that  
\beqq\label{convergence of score}
\call (\, \Delta_n^{(\vth)} \mid P^\vth  )  \;\lra\;  \caln (\, 0 \,,\, \cali^{(\vth)}  \,) \quad\mbox{(weak convergence in $\bbr$, $\nto$)} \;. 
\eeqq

\vskip0.8cm
In the special case $\si(\cdot)\equiv 1$, the Fisher information (\ref{fisherinformation}) coincides with [G 88, p.\ 289] and [IH 81, p.\ 209]. Theorem 1.1 will be proved at the end of section 2 below. 

Under 'discontinuous signal' hypotheses, the rate is $n^{-2}$ (thus essentially faster than the rate $n^{-1}$ in [K 04, section 3.4] for time homogenous ergodic diffusions where the drift has jumps at parameter dependent positions, and also essentially faster than the rate $n^{-1}$ in [HK 10] for time inhomogenous periodic settings with known periodicity where a signal in the drift has jumps at parameter dependent positions), and the nature of the limit experiment, see above, is very different from the Gaussian shift in theorem 1.1. Below, $( \wt W (u) )_{u\in\bbr}$ is double sided standard Brownian motion.\\

{\bf 1.2 Theorem : } Under hypotheses {\bf (H0')}+{\bf (H1')}+{\bf (H2)} consider the sequence of experiments (\ref{sequence_of_experiments}) and local models at $\vth\in\Theta$
\beqq\label{IHlocalmodels}
\left\{ P^{\vth + n^{-2} h}_{n} : h\in\Theta_{\vth,n} \right\}  \quad,\quad n\to\infty 
\eeqq
where  $\Theta_{\vth,n}:= \{ h\in\bbr: \vth + n^{-2} h \in\Theta\}$. For all  $\vth\in\Theta$, we have convergence under $P^\vth$ as $\nto$ of 
$$
\left( L_n^{ (\vth + n^{-2} h) / \vth } \right)_{h\in\Theta_{\vth,n}}
$$
in the sense of finite-dimensional distributions to 
\beqq\label{IHlimitexperiment}
\left( \wt L^{ h / 0 } \right)_{h\in\bbr}  \quad,\quad  L^{ h / 0 } \;:=\; \exp\left\{\, \wt W (\, h\, J(\vth,r,\rho) \,) \;-\; \frac12\, |h|\; J(\vth,r,\rho) \,\right\} 
\eeqq
where the scaling constant  
\beqq\label{new_constants-1}
J(\vth, r,\rho)  \;:=\;  \frac{1}{2\vth^2}\, \sum_{j=1}^\ell\;  \rho_j^2\; [\mu^\vth P^\vth_{0,r_j\vth}](\,\frac{1}{\si^2}\,)  
\eeqq
depends on the collection $r=(r_1,\ldots,r_\ell)$ of jump times and $\rho=(\rho_1,\ldots,\rho_\ell)$ of jump heights in the $1$-periodic signal $S_0(\cdot)$.\\

Theorem 1.2 will be proved at the end of section 3 below.\\

\section*{2. Proofs: laws of large numbers, application to smooth signals}

We start with those parts of the proofs which are common to theorems 1.1 and 1.2; at the end of this section, we will prove theorem 1.1. 

Let the $1$-periodic function $S_0(\cdot)$ be piecewise continuous on $[0,1]$. 
Under $P^\vth$, we have a $\vth$-periodic structure in the semigroup associated to the canonical process $(\eta_t)_{t\ge 0}$  
\beqq\label{semigroup-periodicity}
P^\vth_{s,t}(x,dy) \;=\; P^\vth_{s+k\vth,t+k\vth}(x,dy) \quad\mbox{for $0\le s<t <\infty$, for $k\in\bbz$ such that $s+k\vth\ge 0$}
\eeqq
and thus can define on the canonical path space $\,\left( C([0,\infty)),\calc([0,\infty)) \right)\,$ for solutions of (\ref{process})
\begin{itemize}\item 
a Markov chain $\mathbb{X}^\vth$ of path segments 
$$
\mathbb{X}^\vth =  \left( \mathbb{X}^\vth_k \right)_{k\in\bbn_0}  \quad,\quad   \mathbb{X}^\vth_k := \left( \eta_{(k-1)\vth+s} \right)_{0\le s\le \vth}
$$
which takes values in the space $\left( C([0,\vth]),\calc([0,\vth]) \right)$ of continuous functions $[0,\vth]\to\bbr$  (starting from $\mathbb{X}^\vth_0\equiv\al_0$, for some $\al_0\in C([0,\vth])$ with terminal value $\al_0(\vth)=x_0$); 
\item
a grid chain $X^\vth$ 
$$
X^\vth = (X^\vth_k)_{k\in\bbn_0}  \quad,\quad   X^\vth_k := \eta_{k\vth}  
$$ 
\end{itemize}
which both are time homogeneous Markov chains under $P^\vth$, by (\ref{semigroup-periodicity}). Condition {\bf (H2)} implies as in [HK 10, section 2] that the segment chain $\mathbb{X}^\vth$ is positive Harris recurrent under $P^\vth$. We write $m^\vth$ for the invariant measure on $\left( C([0,\vth]) , \calc([0,\vth]) \right)$ of the segment chain which is determined uniquely from 
\beqq\label{invmeasureTper-1}
\left\{ \begin{array}{l}
\mbox{for arbitrary $0=t_0<t_1<\ldots<t_{\ell}<t_{\ell+1}=\vth$ and $A_i\in\calb(\bbr)$} \;, \\
m^\vth\left( \left\{ \pi_{t_i}\in A_i \,,\,0\le i\le \ell{+}1 \right\}\right) \quad\mbox{is given by}\\
\int\ldots\int \mu^\vth(dx_0)\, 1_{A_0}(x_0)\, \prod_{i=0}^\ell P^\vth_{t_i,t_{i+1}}(x_i,dx_{i+1})\, 1_{A_{i+1}}(x_{i+1}) \;. 
\end{array} \right.
\eeqq
Moreover, as a consequence of  {\bf (H2)}, we have the following strong law of large numbers ([HK 10, theorem 2.1]): for every $\vth\in\Theta$ and every increasing process $A^\vth = (A^\vth_t)_{t\ge 0}$ with the property
\beqq\label{property}
\left\{ \begin{array}{l}
\mbox{there is some function $F^\vth:C([0,\vth])\to\bbr$, nonnegative, $\calc([0,\vth])$-measurable, }\\
\mbox{satisfying}\quad   m^\vth(F^\vth) := \int_{C([0,\vth])} F^\vth(\al)\, m^\vth(d\al) \;<\;  \infty \;, \;\;\mbox{such that} \\
A^\vth_{k\vth} \;=\; \sum\limits_{j=1}^k F^\vth(\mathbb{X}^\vth_j)  
\;=\; \sum\limits_{j=1}^k  F^\vth\left(\, (\eta_{(k-1)\vth+s})_{0\le s\le \vth} \,\right) \;,\; k\ge 1  
\end{array} \right.
\eeqq
we have $P^\vth$-almost surely 
\beqq\label{SLLN}
\lim_{k\to\infty}\, \frac{1}{k}\; \sum\limits_{j=1}^k F^\vth(\mathbb{X}^\vth_j) \;=\; m^\vth(F^\vth) 
\quad\mbox{and}\quad
\lim_{t\to\infty}\, \frac{1}{t}\; A^\vth_t \;=\; \frac{1}{\vth}\; m^\vth(F^\vth)  \;. 
\eeqq
On this basis, we have the following.\\

{\bf 2.1 Lemma: } Assume {\bf (H2)} and $\si(\cdot)$ strictly positive. For any function $f:\bbr\to[0,\infty)$ which is measurable, $1$-periodic, bounded, for every $\vth\in\Theta$, define an increasing process $A^\vth$ by 
\beqq\label{increasingprocess}
A^\vth_t \;\;:=\;\; \int_0^t \frac{f(\frac1\vth\, s)}{\si^2(\eta_s)}\; ds \quad,\quad t\ge 0 \;. 
\eeqq
Then we have $P^\vth$-almost surely as $t\to\infty$ 
$$
\frac1t\; A^\vth_t    \quad\lra\quad  
\int_0^1 dv\; f(v)\; [\mu^\vth P^\vth_{0,v\vth}]( \frac{1}{\si^2} ) \;\;=:\;\; C(\vth,f) 
$$
(where the limit may take the value $+\infty$).

\vskip0.5cm
{\bf Proof: } Fix  $\vth\in\Theta$. Assume first that $\si(\cdot)$ is bounded away from $0$. Define a $\si$-finite measure on $(\bbr,\calb(\bbr))$ 
$$
\Lambda^\vth(ds) \;:=\; f(\frac1\vth\, s)\, ds 
$$
which is $\vth$-periodic in the sense that $\Lambda^\vth(B)=\Lambda^\vth(B+k\vth)$ for  $B\in\calb(\bbr)$ and  $k\in\bbz$, and a functional 
$$
F^\vth : \quad 
C([0,\vth]) \;\ni\;\; \al  \quad\lra\quad   \int_0^\vth  \Lambda^\vth(ds)\, \frac{1}{\si^2(\al(s))}   \;\;<\; \infty \;. 
$$
For $A^\vth$ defined in (\ref{increasingprocess}) we have as in (\ref{property}) 
$$
A^\vth_{k\vth} \;\;=\;\; \sum_{j=1}^k F^\vth(\mathbb{X}^\vth_j) \;\;,\;\; k\ge 1
$$
and apply the strong law of large numbers (\ref{SLLN}). Calculating  
$$
m^\vth(F^\vth)  
\;\;=\;\;   \int_0^\vth  \Lambda^\vth(ds) \int_{C([0,\vth])}m^\vth(d\al)\, (\frac{1}{\si^2(\al(s))}) 
\;\;=\;\;   \int_0^\vth  \Lambda^\vth(ds) \int_{\bbr} [\mu^\vth P^\vth_{0,s}](dy) (\frac{1}{\si^2(y)}) 
$$
the limit of $\frac1t\; A^\vth_t$ under $P^\vth$ as $t\to\infty$ 
$$
\frac1\vth\, m^\vth(F^\vth) \;=\;
\frac{1}{\vth} \int_0^\vth ds\; f(\frac1\vth\, s)\; [\mu^\vth P^\vth_{0,s}]( \frac{1}{\si^2} ) 
$$
equals $C(\vth,f)$ as asserted. The lemma is proved when $\si(\cdot)$ is bounded away from $0$; in the general case we replace $\si$  by $\si\wedge \delta$, $\,\delta>0$, and let $\delta$ tend to $0$. \halmos\\

In the following, we will assume that $\si(\cdot)$ is bounded away from $0$; this guarantees for finite limits $C(\vth,f)$ under arbitrary $\vth$ whenever $f$ is bounded. In the present section, this is merely for convenience, but will be essential in section 3 below (in steps 5) and 6) of the proof of lemma 3.1). \\

{\bf 2.2 Lemma: } Assume {\bf (H1)}+{\bf (H2)}, consider $f:\bbr\to[0,\infty)$  measurable, $1$-periodic, and bounded. Consider a function $H:[0,\infty)\to[0,\infty)$ which is c\`adl\`ag increasing, and which varies regularly at $\infty$ with index $\rho>0$ ([BGT 87]). Then for every $\vth\in\Theta$ 
$$
\frac{1+\rho}{t\, H(t)}\; \int_0^t H(s)\; \frac{f(\frac1\vth\, s)}{\si^2(\eta_s)}\; ds    \quad\lra\quad  
C(\vth,f) \quad\quad\mbox{$P^\vth$-almost surely as $t\to\infty$}  
$$
where the limit is $C(\vth,f)$, as defined in lemma 2.1, is finite.

\vskip0.5cm
{\bf Proof: } We have from [BGT 87, theorem 1.6.4 on p.\ 33]
\beqq\label{integralformula}
\int_{t_0}^t s\, dH(s) \;\;\sim\;\;   \frac{\rho}{1+\rho}\, t\,H(t)   \quad\mbox{as}\;\; t\to\infty    
\eeqq
for arbitrary $0<t_0<\infty$ fixed. Fix $\vth\in\Theta$. The paths of $t\to A^\vth_t$ defined in (\ref{increasingprocess}) being  continuous, Stieltjes product formula gives
$$
\int_0^t H(s)\, dA^\vth_s \;\;=\;\; H(t)\, A^\vth_t \;-\; \int_0^t A^\vth_s\, dH(s) 
$$
for all $t\ge 0$. Both terms on the r.h.s.\ are of order $t H(t)$, as a consequence of $P^\vth$-almost sure convergence of $\frac1t A^\vth_t$ as $t\to\infty$ according to lemma 2.1: for the first term, this is obvious since
$$
A^\vth_t \,H(t)  \quad\sim\quad  C(\vth,f)\;  t\, H(t)  
\quad\mbox{$P^\vth$-almost surely as $t\to\infty$}   
$$
by lemma 2.1; for the second term, we deduce from (\ref{integralformula}) and lemma 2.1
$$
\int_0^t A^\vth_s\, dH(s)   \;\;\sim\;\;  C(\vth,f) \int_{t_0}^t s\, dH(s) \;\;\sim\;\; C(\vth,f)\; \frac{\rho}{1+\rho}\, t\, H(t)  
\quad\mbox{$P^\vth$-almost surely as $t\to\infty$} \;. 
$$
Taking differences, the assertion follows. \halmos\\

Write $\bbd$ for the Skorohod path space of c\`adl\`ag functions $[0,\infty)\to\bbr$. Under assumption {\bf (H0)} we write $\dot S(\vth,t)$, $\ddot S(\vth,t)$ for the derivatives of $S(\vth,t)$ with respect to the parameter $\vth$. \\

{\bf 2.3 Lemma: } Assume {\bf (H0)}+{\bf (H1)}+{\bf (H2)}. Then for all $\vth\in\Theta$, we have weak convergence in $\bbd$ of 
$$
M^n  \;:=\;  \left( \sqrt{\frac{3}{n^3}} \int_0^{tn} \frac{\dot S(\vth,s)}{\si(\eta_s)}\; dB^\vth_s\, \right)_{t\ge 0} 
\quad,\quad \nto
$$
under $P^\vth$ to Brownian motion $B \circ \Phi^\vth$ time changed by $\Phi^\vth$ 
$$
t  \quad\lra\quad  \Phi^\vth(t) \;:=\;  \frac{t^3}{\vth^4}\, C(\vth,[S'_0]^2) 
$$
where  $C(\vth,[S'_0]^2)$ is the limit defined in lemma 2.1 for $f:=[S'_0]^2$. 

\vskip0.5cm
{\bf Proof: } The parametrization (\ref{parametrization}) gives 
$\,\dot S(\vth,t) = \frac{d}{d\vth} S_0(\,\frac1\vth t\,) = -\,\frac{t}{\vth^2}\, S'_0(\,\frac1\vth t\,)\,$.  
Applying lemma 2.2 with $H(t)=t^2$ and $f=[S'_0]^2$ under $P^\vth$ yields 
\beqq\label{proofstep}
\langle M^n \rangle_t \;=\; 
\frac{3}{n^3} \int_0^{tn} \frac{[\dot S(\vth,s)]^2}{\si^2(\eta_s)}\, ds  \;=\; 
\frac{1}{\vth^4}\, \frac{3}{n^3}\int_0^{tn} s^2\; \frac{ [S_0']^2(\,\frac1\vth s\,)}{\si^2(\eta_s)}\, ds 
\quad\lra\quad  \frac{t^3}{\vth^4}\, C(\vth,[S'_0]^2)
\eeqq
$P^\vth$-almost surely as $n\to\infty$, for every $t>0$ fixed. Then the martingale convergence theorem  (cf.\ Jacod and Shiryaev [JS 87, VII.3.22]) applies and gives the result. \halmos\\

{\bf 2.4 Lemma: } Under {\bf (H0)}+{\bf (H1)}+{\bf (H2)}, consider sequences in the parameter space
$$
\vth_n \;:=\;  \vth \;+\; \sqrt{\frac{3}{n^3}}\; h_n \quad,\quad n\ge n_0
$$
defined with respect to some bounded sequence $(h_n)_n$ in $\bbr$ and some fixed reference point $\vth\in\Theta$. Then as $\nto$, the following processes under $P^\vth$
\beam 
N^{n,h_n}_t  &:=&  \int_0^{tn} \frac{ S(\vth_n,s) - S(\vth,s) - (\vth_n-\vth)\, \dot S(\vth,s) }{ \si(\eta_s) }\; dB^\vth_s 
\quad,\quad t\ge 0  \label{reste-1} \\ 
U^{n,h_n}_t  &:=&  \int_0^{tn}  \frac{ [\, S(\vth_n,s) - S(\vth,s) - (\vth_n-\vth)\, \dot S(\vth,s) \,]^2 }{ \si^2(\eta_s) }    ds  
\quad,\quad t\ge 0   \label{reste-2}  \\
V^{n,h_n}_t  &:=&  \int_0^{tn} \frac{ [S(\vth_n,s) - S(\vth,s) - (\vth_n-\vth)\, \dot S(\vth,s)][(\vth_n-\vth)\, \dot S(\vth,s)] }{ \si^2(\eta_s) }\; ds 
\quad,\quad t\ge 0   \label{reste-3}
\eeam
vanish uniformly over compact $t$-intervals as $\nto$.

\vskip0.5cm
{\bf Proof: } Fix $\vth\in\Theta$, write $\,\delta_n=\sqrt{\frac{3}{n^3}}\,$,  and consider sequences $(h_n)_n$ which are bounded by $c$. 

1) Exploiting {\bf (H0)} and (\ref{parametrization}), we show that there are constants $c(\vth,c)$ such that  
\beqq\label{difference-1}
\left| S(\vth_n,s) - S(\vth,s) - (\vth_n-\vth)\, \dot S(\vth,s) \right| \;\;\le\;\; c(\vth,c)\; \frac{3}{n^3}\; (1+s^2) 
\quad,\quad s\ge 0
\eeqq
for all $n$ large enough. We write 
$$
S(\vth_n,s) - S(\vth,s) - (\vth_n-\vth)\, \dot S(\vth,s) \;=\; \frac12\, (\vth_n-\vth)^2\, \ddot S(\, \zeta_{\vth,n,h_n}(s) \,)
$$
with suitable $\zeta_{\vth,n,h_n}(s)$ between $\vth$ and $\vth_n\to \vth$. From  
$\ddot S(\zeta,s) = S_0''(\, \frac1\zeta\, s \,)\frac{s^2}{\zeta^4} + S_0'(\, \frac1\zeta\, s \,)\frac{2\, s}{\zeta^3}$ 
and some upper bound $K$ for $|2\,S_0'|$ and $|S_0''|$,  the l.h.s.\ of (\ref{difference-1}) is smaller than  
$$
\frac{K\, c^2}{2}\; \frac{3}{n^3}\,\left\{ \left(\frac{ s^2 }{ [\vth-\delta_{n_0}c]^4 }+\frac{ s }{ [\vth-\delta_{n_0}c]^3 }\right)\right\}  
$$
uniformly in $n\ge n_0$. Together with $s\vee s^2\le 1+s^2$ on $s\ge 0$, this gives (\ref{difference-1}). 

2) We consider the process $U^{n,h_n}$ in (\ref{reste-2}).  For $t>0$ fixed, the bound (\ref{difference-1}) shows
$$
U^{n,h_n}_t  \;\;\le\;\; c^2(\vth,c)\;\frac{9}{n^6}\; \int_0^{tn} \frac{ (1+s^2)^2 }{ \si^2(\eta_s) }\; ds  
$$
and the assertion is immediate from lemma 2.2 with $H(s)=(1+s^2)^2$ and $f\equiv 1$. 

3) For $t>0$ fixed, we consider angle brackets $\left\langle N^{n,h_n} \right\rangle = U^{n,h_n}$ under $P^\vth$    
and make use of step 2). Burkholder-Davis-Gundy inequality ([IW 89, p.\ 110]) then shows that the $P^\vth$-martingales  $N^{n,h_n}$ vanish uniformly over compact $t$-intervals as $\nto$ under $P^\vth$. 

4) Finally, processes (\ref{reste-3}) vanish uniformly over compact $t$-intervals as $\nto$ under $P^\vth$, by Cauchy-Schwarz inequality combined with (\ref{proofstep}) and step 2). \halmos\\

{\bf 2.5 Proof of theorem 1.1: } Fix $\vth\in\Theta$. 
We start with local scale $\,\delta_n=\sqrt{\frac{3}{n^3}\,}$ instead of $\sqrt{\frac{1}{n^3}\,}$. Then with all notations as in lemma 2.4 except that we write $\wt h$ for the local parameter, for bounded sequences $(\wt h_n)_n$ in $\bbr$ 
$$
\log L_{nt}^{\vth_n / \vth } \;=\; \int_0^{nt} \frac{S(\vth_n,s)-S(\vth,s)}{\si(\eta_s)}\, dB^{(\vth)}_s 
\;-\; \frac12 \int_0^{nt} \left( \frac{S(\vth_n,s)-S(\vth,s)}{\si(\eta_s)} \right)^2 ds 
$$
under $P^\vth$. Adding $\pm (\vth_n-\vth) \dot S(\vth,\cdot)$ in the numerators of the integrands,  we separate leading terms
\beqq\label{leadingterm0}
\wt h_n\, \wt \Delta^n_t(\vth)  \;-\; \frac12 \wt h_n^2\, \wt I^n_t(\vth) 
\eeqq
defined by  
\beqq\label{leadingterm1}
\wt \Delta^n_t(\vth) \;:=\;  \sqrt{\frac{3}{n^3}} \int_0^{tn} \frac{\dot S(\vth,s)}{\si(\eta_s)}\, dB^{(\vth)}_s   \;\;,\;\;   t\ge 0 
\eeqq
\beqq\label{leadingterm2}
\wt I^n_t(\vth) \;:=\; \frac{3}{n^3} \int_0^{tn} \frac{ [\, \dot S(\vth,s) \,]^2 }{\si^2(\eta_s)}\; ds  \;\;,\;\;   t\ge 0
\eeqq
from  remainder terms as defined in (\ref{reste-1})--(\ref{reste-3})
$$
N^{n,\wt h_n}_t \;-\; \frac12  U^{n,\wt h_n}_t   \;-\; V^{n,\wt h_n}_t 
$$
which vanish under $P^\vth$ as $\nto$, uniformly on compact $t$-intervals, by lemma 2.4. Now (\ref{proofstep}) allows to replace $\wt I^n_t(\vth)$ by its limit $\frac{t^3}{\vth^4}\, C(\vth,[S_0']^2)=:\Phi^\vth(t)$ under $P^\vth$ as $\nto$, and lemma 2.3 gives weak convergence of the martingales $\wt \Delta^\vth$ in (\ref{leadingterm0}) under $P^\vth$ to $B\circ\Phi^\vth$ as $\nto$. 
We thus have proved a quadratic expansion of log-likelihood ratios with local parameter $\wt h\in\bbr$ and local scale $\sqrt{\frac{3}{n^3}\,}$. Viewing  $\sqrt{3}\, \wt h =: h$ as local parameter, we get (\ref{quadraticdecomposition})+(\ref{score})+(\ref{convergence of score}) as stated in theorem~1.1.  \halmos \\

\section*{3. Proofs: discontinuous signals}

We continue with 'general' signals $S_0(\cdot)$ as in the beginning of section 2. At the end of the present section, we will prove theorem 1.2.\\

{\bf 3.1 Lemma: } Assume {\bf (H1')}+{\bf (H2)}. Consider sequences $m \,\sim\, c\,n$ for suitable $0<c<\infty$. Fix $0<r<1$ and $h\in\bbr$. Then under $P^\vth$ as $\nto$ 
\beam\label{approximation-1}
\sum_{k=0}^m \int_{(\vth+h/n^2)(k+r)}^{\vth(k+r)} \frac{1}{\si^2(\eta_s)} ds 
&=& |h|\, \frac{1}{n^2} \sum_{k=0}^m\; k\; \frac{1}{\si^2(\eta_{\vth(k+r)})} \;+\; o_{P^\vth}(1) \quad\quad\mbox{if $h<0$} \;, 
\\ \label{approximation-2}
\sum_{k=0}^m \int_{\vth(k+r)}^{(\vth+h/n^2)(k+r)} \frac{1}{\si^2(\eta_s)} ds 
&=& h\; \frac{1}{n^2} \sum_{k=0}^m\; k\; \frac{1}{\si^2(\eta_{\vth(k+r)})} \;+\; o_{P^\vth}(1) \quad\quad\mbox{if $h>0$}  \;.  
\eeam
The leading term on the r.h.s.\ of (\ref{approximation-1})+(\ref{approximation-2}) has the property 
\beqq\label{approximation-4}
\frac1m \sum_{k=0}^m \; \frac{1}{\si^2(\eta_{\vth(k+r)})}  \quad\lra\quad  [\mu^\vth P^\vth_{0,r\vth}](\frac{1}{\si^2})   
\quad\quad\mbox{$P^\vth$-almost surely as $\nto$}
\eeqq
and for increasing sequences $(H(n))_n$ which vary regularly at $\infty$ with index $\rho>0$ we have 
\beqq\label{approximation-3}
\frac{1+\rho}{m\, H(m)} \sum_{k=0}^m \;H(k)\; \frac{1}{\si^2(\eta_{\vth(k+r)})}  \quad\lra\quad  [\mu^\vth P^\vth_{0,r\vth}](\frac{1}{\si^2})   \quad\quad\mbox{$P^\vth$-a.s.\ as $\nto$} \;. 
\eeqq

\vskip0.5cm
{\bf Proof: } 1) Fix $\vth\in\Theta$. Imitating the proof of lemma 2.1 with 
$$
A^\vth_t \;:=\; \int_0^t \frac{1}{\si^2(\eta_s)}\, \Lambda^\vth(ds)   \quad,\quad  \Lambda^\vth := \sum\limits_{k\in\bbz}\epsilon_{\vth(k+r)} \quad,\quad  F^\vth\left(\, (\al(s))_{0\le s\le \vth} \,\right) := \frac{1}{\si^2(\al(r\vth))} 
$$ 
under $P^\vth$ we obtain (\ref{approximation-4}) from the strong law of large numbers (\ref{SLLN}). To prove (\ref{approximation-3}), we embed the sequence $(H(n))_n$ into a c\`adl\`ag increasing function $\wt H:[0,\infty)\to[0,\infty)$ such that  $\wt H( \vth(k+r) )=H(k)$ for all $k$. Then $\wt H(\cdot)$ is regularly varying at $\infty$ with index $\rho$.  Based on (\ref{approximation-4}) and on the Stieltjes product formula for c\`adl\`ag paths $t\to A_t^\vth$ 
$$
\int_0^t \wt H(s)\, dA^\vth_s \;\;=\;\; \wt H(t)\, A^\vth_t \;-\; \int_0^t A^\vth_{s-}\, d\wt H(s)  
$$
we imitate the proof of lemma 2.2 to obtain 
$$
\frac{1+\rho}{t\;\wt H(t)}\; \int_0^t \wt H (s)\, dA^\vth(s)   \quad\lra\quad   \frac1\vth\; m^\vth(F^\vth) = 
\frac1\vth\; [\mu^\vth P^\vth_{0,r\vth}](\, \frac{1}{\si^2} \,)
$$  
$P^\vth$-almost surely as $t\to\infty$. Cancelling a factor $\frac1\vth$ on both sides gives (\ref{approximation-3}).

2) We shall apply the following exponential inequality (\ref{exponential_inequality}), due to [B 05] and adapted to time inhomogeneous diffusions in [HK 10, lemma 3.1]. With $\vth$-periodic drift as in (\ref{process}), this inequality requires Lipschitz conditions on $b(\cdot)$ and $\si(\cdot)$, boundedness of $S_0(\cdot)$, and boundedness of $\si(\cdot)$ which is guaranteed by {\bf (H1')}.  
Fix $0<\la<\frac12$ and $\frac12<\kappa<1{-}\lambda$.  Then for the canonical process $\eta$ under $P^\vth$, there is some $\Delta_0>0$ such that  
\beqq\label{exponential_inequality}
P^\vth\left(\, \sup_{t_1\le t\le t_1+\Delta}|\eta_t-\eta_{t_1}|>\Delta^\lambda \;,\; |\eta_{t_1}|\le \left(\frac1\Delta\right)^\kappa\,\right)  \;\;\le\;\; c_1\,\cdot\, \exp\left\{ \,-\, c_2 \left(\frac1\Delta\right)^{1-2\lambda}\right\}  
\eeqq
holds for all $0\le t_1<\infty$ and all $0<\Delta<\Delta_0$, with positive constants $c_1$ and $c_2$ which do not depend on $t_1\ge 0$ or on $\Delta\in(0,\Delta_0)$.

3) Below we give a detailed proof for assertion (\ref{approximation-1}) which corresponds to the case $h<0$ (the proof of (\ref{approximation-2}) for $h>0$ is then similiar, and slightly simpler). Integration on the l.h.s.\ in (\ref{approximation-1}) is over intervals 
\beqq\label{collection_of_intervals}
[ (\vth+h/n^2)(k+r) ,\vth(k+r) ] \quad,\quad 0\le k\le m  
\eeqq
where $h<0$, $0<r<1$, and $m \sim cn$ as $\nto$ for suitable $0<c<\infty$. Put $\Delta_n:= d\frac{m+1}{n^2}$. Uniformly in $|h|\le d$, for $n$ large enough, the intervals in (\ref{collection_of_intervals}) of length $|h|\frac{k+r}{n^2}$ can be embedded into intervals of equal length 
\beqq\label{larger_intervals}
J_{n,k} := [ \vth(k+r)-\Delta_n , \vth(k+r) ] \quad\subset\quad [\vth k',\vth(k'+2)]   \quad,\quad |k'-k|\le 1  \;. 
\eeqq 
Here the intervals on the r.h.s. of (\ref{larger_intervals}) correspond to twice the periodicity $\vth$ of the canonical process $\eta$ under $P^\vth$. Under {\bf (H2)}, we have Harris recurrence of the $\vth$-segment chain $\mathbb{X}^\vth=(\mathbb{X}^\vth_k)_k$ with invariant probability $m^\vth$ on $C([0,\vth])$. As a consequence, we also have Harris recurrence of a 'bivariate' chain 
$$
\mathbb{X}^{\vth,2}=(\mathbb{X}^{\vth,2}_k)_k  \quad,\quad   \mathbb{X}^{\vth,2}_k := \left( (\eta_{\vth(k+s)} )_{0\le s\le 2} \right)_k
$$
formed by two successive $\vth$-segments; the invariant probability $m^{\vth,2}$ on $C([0,2\vth])$ for the chain $\mathbb{X}^{\vth,2}$ is easily determined, in analogy to (\ref{invmeasureTper-1}) for $\mathbb{X}^\vth$.

4) We imitate step 1) of the proof of theorem 4.1 in [HK 10] to show the following: for arbitrary constants $K>\infty$ and powers $\kappa>0$ 
\beqq\label{auxiliary-1}
\frac{1}{m} \sum_{k=0}^{m} 1_{\{ \sup\limits_{0\le s\le 2} |\eta_{\vth(k+s)}| \,>\, K\, m^\kappa \}} \;\;\lra\;\; 0 
\eeqq
$P^\vth$-almost surely as $\nto$. This follows from the strong law of large numbers for additive functionals of the chain $\mathbb{X}^{\vth,2}$: using  $\,F^\vth: \al \to  1_{ \{\sup\limits_{0\le s\le 2\vth}|\al(s)|\, >\, c\} }\,$ which is defined on $C([0,2\vth])$, we have 
$$
\frac1m \sum_{k=0}^m F^\vth (\, \mathbb{X}^{\vth,2}_k\,)  \quad\lra\quad 
m^{\vth,2} \left\{ \al\in C([0,2\vth]) : \sup\limits_{0\le s\le 2\vth}|\al(s)| > c \right\} 
$$
$P^\vth$-almost surely as $m\to\infty$ where the r.h.s.\ can be made arbitrarily small for large choice of  $c$.

5)  We shall show  
\beqq\label{auxiliary-2}
\frac{1}{m} \sum_{k=0}^m \frac{1}{\si^2(\eta_{(\vth+h/n^2)(k+r)})}    
\;\;=\;\;  \frac{1}{m} \sum_{k=0}^m \frac{1}{\si^2(\eta_{\vth(k+r)})} \;+\; o_{P^\vth}(1) 
\quad\quad\mbox{under $P^\vth$ as $\nto$}
\eeqq
for $h<0$, $0<r<1$, and $m \sim c n$ as $\nto$. 

The idea is to control fluctuations of the canonical process $\eta$ over intervals of identical length $\Delta_n = d\frac{m+1}{n^2}=O(\frac1n)$ thanks to the exponential inequality (\ref{exponential_inequality}), similiar to step 2) of the proof of theorem 4.1 of [HK 10]. For $|h|\le d$ and $n$ large enough, with notations of (\ref{larger_intervals}), we embed  
$$
[ (\vth+h/n^2)(k+r) ,\vth(k+r) ]  \quad\subset\quad  [ \vth(k+r)-\Delta_n , \vth(k+r) ] = J_{n,k} 
\quad,\quad 0\le k\le m   
$$
and consider the $k$-th summand contributing to the difference in (\ref{auxiliary-2})  
\beqq\label{auxiliary-3}
\left| \frac{1}{\si^2(\eta_{\vth(k+r)})} -  \frac{1}{\si^2(\eta_{(\vth+h/n^2)(k+r)})}  \right|  \;. 
\eeqq
This summand admits --~since $\si(\cdot)$ is Lipschitz and bounded away from $0$ and $\infty$~--  bounds of type  
\beao
&&d_1\cdot 1_{ \{\, \sup\limits_{0\le s\le 2\vth}\, |\eta_{\vth(k-1)+s}| \;\;>\;\;   \left(\frac{1}{\Delta_n}\right)^\kappa  \,\}  }      \\
&&+\quad d_1\cdot 1_{ \{\, \sup\limits_{ t\in J_{n,k} }\, |\eta_t-\eta_{\vth(k+r)-\Delta_n }|  \;\;>\;\; \Delta_n^\la  \; \;,\; \; |\eta_{ \vth(k+r)-\Delta_n}|  \;\;\le\;\;   \left(\frac{1}{\Delta_n}\right)^\kappa     \,\}  }  \\
 &&+\quad d_2\, \Delta_n^\la    \cdot    1_{ \{\,  \sup\limits_{ t\in J_{n,k} }\, |\eta_t-\eta_{\vth(k+r)-\Delta_n }|  \;\; \le\;\; \Delta_n^\la  \,\}  } 
\eeao 
with suitable constants $d_1$, $d_2$; here $d_2$ involves the Lipschitz constant for $\si(\cdot)$, and $d_1$ the lower bound for $\si(\cdot)$. By the type of bound in the first line combined with $m\sim cn$ and step 4), we see that 
$$
 \frac{1}{m} \sum_{j=0}^{m-1}  \left| \frac{1}{\si^2(\eta_{\vth(k+r)})} -  \frac{1}{\si^2(\eta_{(\vth+h/n^2)(k+r)})}  
\right| \; 1_{ \{\, \sup\limits_{0\le s\le 2\vth}\, |\eta_{\vth(k-1)+s}| \;\;>\;\;   \left(\frac{1}{\Delta_n}\right)^\kappa  \,\}  }   
$$
 vanishes almost surely under $P^\vth$ as $\nto$. Next, the exponential inequality in step 2) on the intervals $J_{n,k}$ of length $\Delta_n$ (i.e.\ with $t_1 = \vth(k+r)-\Delta_n$ for $k=0,1,\ldots, m$) shows that  
 $$
 P^\vth \left(\,    \sup\limits_{ t\in J_{n,k}}\, |\eta_t-\eta_{ \vth(k+r)-\Delta_n}|  \;> \; \Delta_n^\la  \; \;,\; \; |\eta_{ \vth(k+r)-\Delta_n}|  \;\le \;   \left(\frac{1}{\Delta_n}\right)^\kappa  \; ,\;\mbox{some}\; \; k=0,1,\ldots,m \,\right)   
 $$
 vanishes under $P^\vth$ as $\nto$. Hence, by the type of bound in the second line,  the probability under $\vth$ to find any strictly positive summand in the sum 
 $$
 \frac{1}{m} \sum_{j=0}^{m-1}  \left| \frac{1}{\si^2(\eta_{\vth(k+r)})} -  \frac{1}{\si^2(\eta_{(\vth+h/n^2)(k+r)})}  
\right| \;  1_{ \{\, \sup\limits_{ t\in J_{n,k} }\, |\eta_t-\eta_{\vth(k+r)-\Delta_n }|  \;\;>\;\; \Delta_n^\la  \; \;,\; \; |\eta_{ \vth(k+r)-\Delta_n}|  \;\;\le\;\;   \left(\frac{1}{\Delta_n}\right)^\kappa     \,\}  }  
$$
tends to $0$ as $\nto$: hence the sum itself vanishes under $P^\vth$ as $\nto$.  Finally, by the type of bound in the third line, we are left to consider averages 
$$
\frac{1}{m} \sum_{j=0}^{m-1}  \left| \frac{1}{\si^2(\eta_{\vth(k+r)})} -  \frac{1}{\si^2(\eta_{(\vth+h/n^2)(k+r)})}  
\right| \;  1_{ \{\, \sup\limits_{ t\in J_{n,k} }\, |\eta_t-\eta_{\vth(k+r)-\Delta_n }|  \;\;\le\;\; \Delta_n^\la      \,\}  }  
$$
which are bounded by $d_2 \Delta_n^\la$, and thus vanish  as $\nto$. We have proved (\ref{auxiliary-2}).

6) Next we prove that  under $P^\vth$, as $\nto$,   
\beqq\label{auxiliary-4}
\sum_{k=0}^m \int_{(\vth+h/n^2)(k+r)}^{\vth(k+r)} \frac{1}{\si^2(\eta_s)}\, ds   
\;\;=\;\; |h|\, \frac{1}{n^2} \sum_{k=0}^m  \;(k+r)\,\frac{1}{\si^2(\eta_{(\vth+h/n^2)(k+r)})} \;+\; o_{P^\vth}(1)
\eeqq 
for $h<0$, $0<r<1$, and $m \sim c n$ as $\nto$. 

For $k=0,1,\ldots, m$, consider in the difference between l.h.s.\ and r.h.s.\ of (\ref{auxiliary-4}) a $k$-th summand 
$$
\int_{(\vth+h/n^2)(k+r)}^{\vth(k+r)} \left| \frac{1}{\si^2(\eta_{s})} -  \frac{1}{\si^2(\eta_{(\vth+h/n^2)(k+r)})}  \right| ds \;. 
$$
Exploiting $m\sim cn$ and thus $\Delta_n=O(\frac1m)$ as $\nto$, we embed the interval of integration into the larger $J_{n,k}$ of (\ref{larger_intervals}), uniformly in $|h|\le d$, and have for the $k$-th summand bounds of type  
\beao
&&\frac1m\, d_1\cdot 1_{ \{\, \sup\limits_{0\le s\le 2\vth}\, |\eta_{\vth(k-1)+s}| \;\;>\;\;   \left(\frac{1}{\Delta_n}\right)^\kappa  \,\}  }      \\
&&+\quad \frac1m\, d_1\cdot 1_{ \{\, \sup\limits_{ t\in J_{n,k} }\, |\eta_t-\eta_{\vth(k+r)-\Delta_n }|  \;\;>\;\; \Delta_n^\la  \; \;,\; \; |\eta_{ \vth(k+r)-\Delta_n}|  \;\;\le\;\;   \left(\frac{1}{\Delta_n}\right)^\kappa     \,\}  }  \\
 &&+\quad \frac1m\, d_2\, \Delta_n^\la    \cdot    1_{ \{\,  \sup\limits_{ t\in J_{n,k} }\, |\eta_t-\eta_{\vth(k+r)-\Delta_n }|  \;\; \le\;\; \Delta_n^\la  \,\}  } \;. 
\eeao 
The factor $\frac1m$ here in combination with $\sum_{k=0}^m$ in (\ref{auxiliary-4}) allows to proceed in exact analogy to step 5) above to establish (\ref{auxiliary-4}).

7) To finish the proof, we deduce (\ref{approximation-1}) from (\ref{auxiliary-4}). With $H(k):=k+r$, the difference  
\beqq\label{difference-final}
\frac{1}{m\, H(m)}  \left|\, \sum_{k=0}^m  \;H(k)\,\frac{1}{\si^2(\eta_{(\vth+h/n^2)(k+r)})} 
\;-\; \sum_{k=0}^m  \;H(k)\,\frac{1}{\si^2(\eta_{\vth(k+r)})} \,\right|
\eeqq
is trivially smaller than   
$$
\frac1m \sum_{k=0}^m \left|\, \frac{1}{\si^2(\eta_{(\vth+h/n^2)(k+r)})}  \;-\; \frac{1}{\si^2(\eta_{\vth(k+r)})} \,\right|
$$
which vanishes $P^\vth$-almost surely as $\nto$ by (\ref{auxiliary-2}). Exploiting $m\sim c n$ and (\ref{approximation-3}), we replace the factor in front of the  difference (\ref{difference-final}) by $\frac{1}{n^2}$, and see --using (\ref{approximation-4}) for the last comparison -- that 
$$
\frac{1}{n^2}  \sum_{k=0}^m  \;(k+r)\,\frac{1}{\si^2(\eta_{(\vth+h/n^2)(k+r)})} 
\quad,\quad  \frac{1}{n^2} \sum_{k=0}^m  \;(k+r)\,\frac{1}{\si^2(\eta_{\vth(k+r)})}  
\quad,\quad  \frac{1}{n^2} \sum_{k=0}^m  \;k\,\frac{1}{\si^2(\eta_{\vth(k+r)})} 
$$
are equivalent under $P^\vth$ as $\nto$. Hence the leading terms on the r.h.s.\ of (\ref{approximation-1}) and (\ref{auxiliary-4}) are equivalent under $P^\vth$ as $\nto$. This etablishes (\ref{approximation-1}) and concludes the proof. \halmos\\

Fix a collection of points $0<r_1<\ldots<r_\ell<1$ and a collection of step functions 
$$
N_j(t) \;:=\;  \sum_{k=0}^\infty 1_{[k+r_j,\infty)}(s) \quad,\quad 1\le j\le \ell    
$$
without common jumps. With respect to some fixed reference point $\vth\in\Theta$  we define for $s\ge 0$
\beqq\label{definition_jotthetahaen}
j^{h,n}_\vth(s) \;:=\; 
\left\{\begin{array}{ll}
\sum\limits_{j=1}^\ell \rho_j \left[ N_j(\,\frac{1}{(\vth+h/n^2)}\,s\,) - N_j(\,\frac{1}{\vth}\,s\,) \right] &  \mbox{in case $h<0$} \\
\sum\limits_{j=1}^\ell \rho_j \left[ N_j(\,\frac{1}{\vth}\,s\,) - N_j(\,\frac{1}{(\vth+h/n^2)}\,s\,) \right] &  \mbox{in case $h>0$} 
\end{array}\right.
\eeqq
for arbitrary real $\rho_j$, $1\le j\le \ell$, and arbitrary $h\in\bbr$, and write as in (\ref{new_constants-1}) 
$$
J(\vth, r,\rho)  \;:=\;  \frac{1}{2\vth^2}\, \sum_{j=1}^\ell\;  \rho_j^2\; [\mu^\vth P^\vth_{0,r_j\vth}](\,\frac{1}{\si^2}\,)  
$$
for $\rho=(\rho_1,\ldots,\rho_\ell)$ and $r=(r_1,\ldots,r_\ell)$.\\

{\bf 3.2 Lemma: } Assume {\bf (H1')}+{\bf (H2)}. For all $\vth\in\Theta$, we have convergence in $P^\vth$-probability  
\beao
\int_0^{tn} j^{h_1,n}_\vth(s)\;  j^{h_2,n}_\vth(s)\; \frac{1}{\si^2(\eta_s)}\; ds 
&\lra&  (|h_1|\wedge|h_2|)\; t^2\;  J(\vth, r,\rho) \quad\mbox{in case ${\rm sgn}(h_i)={\rm sgn}(h_j)$}  \\ 
\int_0^{tn} j^{h_1,n}_\vth(s)\;  j^{h_2,n}_\vth(s)\; \frac{1}{\si^2(\eta_s)}\; ds 
&\lra&  0 \quad\quad\mbox{in case ${\rm sgn}(h_1)\neq {\rm sgn}(h_2)$}  
\eeao
for every $0<t<\infty$ fixed, as $n\to\infty$. \\

{\bf Proof: } Fix $\vth\in\Theta$, fix some $0<t_0<\infty$, fix $h_1,h_2$ such that $|h_i|<d$, $i=1,2$. We consider $t\le t_0$ and define $\Delta_n$ by (\ref{larger_intervals}). 

1) Consider first the case $h_1<h_2<0$. We can choose $n$ large enough to make sure that for all $j=1,\ldots,\ell$, differences of counting functions  
$$
s \quad\lra\quad  N_j(\,\frac{1}{(\vth+h_i/n^2)}\,s\,) - N_j(\,\frac{1}{\vth}\,s\,)  
$$
are $\{0,1\}$-valued in restriction to $[0,t_0n]$, and supported --in restriction to $[0,t_0n]$-- by collections of intervals of a form (\ref{collection_of_intervals})  
\beqq\label{meine_intervalle}
[ (\vth+h_i/n^2)(k+r_j) ,\vth(k+r_j) ] \quad,\quad 0\le k\le m \;\;\mbox{for suitable $m$ with $m \,\sim\, \frac{t_0}{\vth}\, n$}   
\eeqq
with the following property:  whenever $j\neq j'$, no intersections occur between a collection (\ref{meine_intervalle}) of intervals corresponding to $\,r_j\,$ and a collection  (\ref{meine_intervalle}) of intervals corresponding to $\,r_{j'}\,$, irrespectively of the values $|h_i|\le d$ under consideration, provided $n$ is larger than some $n_0(\vth,d,t_0)$. 
In case  $h_1>h_2>0$, we use instead of (\ref{meine_intervalle})
\beqq\label{meine_intervalle_bis}
[ \vth(k+r_j),  (\vth+h_i/n^2)(k+r_j) ] \quad,\quad 0\le k\le m \;\;\mbox{for suitable $m$ with $m \,\sim\, \frac{t_0}{\vth}\, n$}     
\eeqq
and obtain the same conclusion.  

2) Consider $h_1,h_2$ such that $|h_i|\le d$, $i=1,2$, and $n\ge n_0(\vth,d,t_0)$. For $r_j$ in $r=(r_1,\ldots,r_\ell)$ fixed, the intervals in (\ref{meine_intervalle}) or (\ref{meine_intervalle_bis}) induced by $(h_2,r_j)$ are subsets of the corresponding intervals in (\ref{meine_intervalle}) or (\ref{meine_intervalle_bis}) induced by $(h_1,r_j)$ in the two cases $h_1<h_2<0$ and $0<h_2<h_1$, and  have void  intersection in  cases $h_1<0<h_2$ or $h_2<0<h_1$. In virtue of step 1) this implies  
\beqq\label{orthogonalitaetsrelation-0}
j^{h_1,n} \cdot j^{h_2,n} \;\equiv\; 0 \quad\mbox{on $\,[0, t_0 n]\,$ whenever $\,{\rm sgn}(h_1)\neq {\rm sgn}(h_2)$}
\eeqq
in restriction to $[0, t_0 n]$, and
\beqq\label{orthogonalitaetsrelation-1}
j^{h_1,n} \cdot j^{h_2,n} \;=\; \left(j^{ \wt h,n}\right)^2 \;\;\mbox{on $\,[0, t_0 n]\,$ where}\;
\left\{ \begin{array}{ll}
\wt h \;:=\; -(|h_1|\wedge|h_2|) & \mbox{if $h_1<0$, $h_2<0$} \\ 
\wt h \;:=\;  h_1\wedge h_2  & \mbox{if $h_1>0$, $h_2>0$} \;.   
\end{array} \right.
\eeqq

3) For $|h_i|\le d$, $t\le t_0$,  and $n\ge n_0(\vth,d,t_0)$, consider first the case $h_1,h_2<0$. With notation of (\ref{orthogonalitaetsrelation-1}), we have 
\beqq\label{schluss-1}
\int_0^{tn} j^{h_1,n}_\vth(s)\;  j^{h_2,n}_\vth(s)\; \frac{1}{\si^2(\eta_s)}\; ds \quad=\quad 
\int_0^{tn} \left(j^{ \wt h,n}(s)\right)^2\; \frac{1}{\si^2(\eta_s)}\; ds    
\eeqq
for $t\le t_0$. For fixed value of $t$, as $\nto$, the support of the integrand on the r.h.s.\ is by step~1) above a system of mutually disjoint intervals 
$$
[ (\vth+\wt h/n^2)(k+r_j) ,\vth(k+r_j) ] \quad,\quad 0\le k\le m \;\;,\;\; m \,\sim\, \frac{t}{\vth}\, n  \;\;,\;\; 1\le j\le \ell \;. 
$$
Hence, by lemma 3.1 with $m \sim \frac{t}{\vth}\, n$ as $\nto$, we can continue equation (\ref{schluss-1}) under $P^\vth$ by 
\beao
&=&  \sum\limits_{j=1}^\ell \;\rho^2_j\; \int_0^{tn}\left[ N_j(\,\frac{1}{(\vth+\wt h/n^2)}\,s\,) - N_j(\,\frac{1}{\vth}\,s\,) \right]^2 \frac{1}{\si^2(\eta_s)}\; ds \\
&=&  \sum\limits_{j=1}^\ell \;\rho^2_j\; \left( \sum_{k=0}^m \int_{(\vth+\wt h/n^2)(k+r_j)}^{\vth(k+r_j)} \frac{1}{\si^2(\eta_s)}\, ds \right) \\ 
&=&  \sum\limits_{j=1}^\ell \;\rho^2_j\; \left( |\wt h|\, \frac{1}{n^2} \sum_{k=0}^m\; k\; \frac{1}{\si^2(\eta_{\vth(k+r_j)})}\right) \;+\; o_{P^\vth}(1) 
\eeao
as $\nto$. Making use of (\ref{approximation-3}) combined with $m \sim \frac{t}{\vth}\, n$, we determine the limit of the leading terms under $P^\vth$ as $\nto$ as
\beqq\label{schluss-2}
|\wt h|\; \frac{t^2}{2\vth^2}\; \sum\limits_{j=1}^\ell \;\rho^2_j\; [\mu^\vth P^\vth_{0,r_j\vth}](\,\frac{1}{\si^2}\,)
\quad=\quad  (|h_1|\wedge|h_2|)\; t^2\; J(\vth,\rho,r)
\eeqq
where $J(\vth,\rho,r)$ is the constant defined in (\ref{new_constants-1}). Hence in case $h_1,h_2<0$, the proof of the lemma is finished. The proof in case $h_1,h_2>0$ is similiar; when ${\rm sgn}(h_1)\neq {\rm sgn}(h_2)$, the proof was already finished with the orthogonality (\ref{orthogonalitaetsrelation-0}). \halmos\\

The following proposition is the key to convergence of local experiments in the 'discontinuous signals' setting, with Ibragimov and Khasminskii's limit experiment (\ref{IHlimitexperiment}), and plays  in this context the same role which 2nd Le Cam lemmata play for convergence to Gaussian shifts. \\

{\bf 3.3 Proposition: } Under hypotheses {\bf (H0')}+{\bf (H1')}+{\bf (H2)}, for local scale  $n^{-2}$, for arbitrary $t_0<\infty$ and for bounded sequences $(h_n)_n$ in $\bbr$, we have a decomposition of log-likelihood ratio processes in local models (\ref{IHlocalmodels}) at $\vth\in\Theta$ as follows:  
\beqq\label{IHtypedecomposition}
\sup_{t\in [0,t_0]} \left|\; \log L_{tn}^{(\vth + n^{-2} h_n) / \vth}  \;-\;  
\left\{ \int_0^{tn} j^{h_n,n}_\vth(s)\, \frac{1}{\si(\eta_s)}\, dB^\vth_s \;-\; \frac12 \int_0^{tn} \left( j^{h_n,n}_\vth(s) \right)^2 \frac{1}{\si^2(\eta_s)}\, ds \right\} \right| 
\eeqq
vanishes in $P^\vth$-probability as $\nto$. Here $j^{h_n,n}_\vth(\cdot)$ are the deterministic functions defined in (\ref{definition_jotthetahaen}) which enjoy the orthogonality properties (\ref{orthogonalitaetsrelation-0})+(\ref{orthogonalitaetsrelation-1}) above for $n$ large enough. \\

{\bf Proof: } 1) Fix $\vth\in\Theta$ and $0<t_0<\infty$. In local models (\ref{IHlocalmodels}) at $\vth$, we start from (\ref{likelihoodratioprocess})+(\ref{parametrization}) and have 
$$
\log L^{(\vth + n^{-2} h_n) / \vth}_{tn} 
\;=\;  \int_0^{tn} \frac{S_0(\,\frac{1}{\vth+ n^{-2} h_n} s\,)-S_0(\,\frac{1}{\vth} s\,)}{\si(\eta_s)}\; dB^{(\vth)}_s \;-\; \frac12 \int_0^{tn}  \frac{ [ S_0(\,\frac{1}{\vth+ n^{-2} h_n} s\,)-S_0(\,\frac{1}{\vth} s\,) ]^2 }{\si^2(\eta_s)}\; ds 
$$
under $P^\vth$. Below, we consider $|h_n|\le d$, $0\le t\le t_0$ and $n\ge n_0(\vth,d,t_0)$, with notations defined in the proof of lemma 3.2. 

2) By assumption {\bf (H0')}, the $1$-periodic function $S_0(\cdot):[0,\infty)\to\bbr$ has finitely many jumps in $(0,1]$, and is Lipschitz between its  jumps. Let $0<r_1<\ldots<r_\ell\le 1$ denote the jump times of $S_0(\cdot)$ in $(0,1]$, and $\rho_1,\ldots,\rho_\ell$ the corresponding jump heights. From $1$-periodic continuation of 
$$
(0,1]\;\ni\; t \;\lra\; \rho_j\left(1_{[r_j,1]}(t) - t\right)  
$$
define  $1$-periodic functions 
$$
\wt S_j( t ) \;\;:=\;\; \rho_j\sum_{k=0}^\infty 1_{[k+r_j,\infty[}(t) \;\;-\;\; \rho_j\, t   \quad,\quad t\ge 0 
$$ 
for every  $1\le j\le \ell$. They have the property $\wt S_j(k)=0$ for all $k\in\bbn_0$ and all $1\le j\le \ell$. Introducing 
$$
\wt S_0( t ) \;\;:=\;\; S_0( t ) \;-\; \left[ \wt S_1( t ) + \ldots + \wt S_\ell( t ) \right]
$$
which is Lipschitz on $[0,\infty)$ and $1$-periodic, we rearrange the decomposition of $S_0(\cdot)$ in form 
\beqq\label{new_form}
S_0( t ) \;\;=\;\; \sum_{j=1}^\ell \rho_j\, N_j(t) \;\;-\;\; S^c_0( t ) \quad,\quad t\ge 0
\eeqq
where $N_j(t) = \sum_{k=0}^\infty 1_{[k+r_j,\infty[}(t)$ and where $S^c_0( \cdot ):[0,\infty)\to\bbr$ collects contributions to   the decomposition of $S_0(\cdot)$ which are Lipschitz. Write  $L$ for the Lipschitz constant of $S^c_0(\cdot)$ on $[0,\infty)$.

3) With notations (\ref{new_form})+(\ref{definition_jotthetahaen}), restarting at step 1) above, we write differences in (\ref{IHtypedecomposition})  
$$
\log L^{(\vth + n^{-2} h_n) / \vth}_{tn} 
\;-\; \left\{  \int_0^{tn} \frac{j^{h_n,n}_\vth(s)}{\si(\eta_s)}\; dB^{(\vth)}_s \;-\; \frac12 \int_0^{tn} \frac{[ j^{h_n,n}_\vth(s) ]^2}{\si^2(\eta_s)}\; ds \right\} \quad,\quad 0\le t\le t_0
$$
as a sum of two processes: a first one where parametrization comes in through a deterministic function which is Lipschitz in $s$
\beqq\label{term-1}
\int_0^{tn} \frac{S^c_0(\,\frac{1}{\vth+ n^{-2} h_n} s\,)-S^c_0(\,\frac{1}{\vth} s\,)}{\si(\eta_s)}\; dB^{(\vth)}_s \;-\; \frac12 \int_0^{tn}  \frac{ [ S^c_0(\,\frac{1}{\vth+ n^{-2} h_n} s\,)-S^c_0(\,\frac{1}{\vth} s\,) ]^2 }{\si^2(\eta_s)}\; ds 
\eeqq
and a second one containing mixed terms 
\beqq\label{term-2}
\int_0^{tn}  \frac{ [ S^c_0(\,\frac{1}{\vth+ n^{-2} h_n} s\,)-S^c_0(\,\frac{1}{\vth} s\,) ] [j^{h_n,n}_\vth(s)] }{\si^2(\eta_s)}\; ds \;. 
\eeqq
i) Since $|h_n|\le d$ and since $S_0^c(\cdot)$ has Lipschitz constant $L$,  angle bracket terms in (\ref{term-1}) are of order  
$$
O(\frac{1}{n^4}) \int_0^{tn} \frac{s^2}{\si^2(\eta_s)}\; ds \quad,\quad 0\le t\le t_0 \;. 
$$
By {\bf (H1')} --or by lemma 3.1-- they vanish uniformly over  $[0,t_0]$ in $P^\vth$-probability as $\nto$. \\
ii) As a consequence of i), by Burkholder-Davis-Gundy inequality, the martingale terms in (\ref{term-1}) vanish uniformly over $[0,t_0]$ in $P^\vth$-probability as $\nto$. \\
iii) We turn to the mixed terms (\ref{term-2}). Since $|h_n|\le d$ we have inclusions, $h_n$ being positive or negative (cf.\ proof of lemma 3.2)  
$$
j^{h_n,n}_\vth(s) \;\le\;  j^{+d,n}_\vth(s) + j^{-d,n}_\vth(s) \quad,\quad s\in[0,t_0n]
$$
and by lemma 3.2 tightness under $P^\vth$ as $\nto$ of 
$$
\int_0^{tn} \frac{[ j^{+d,n}_\vth(s) ]^2 + [ j^{-d,n}_\vth(s) ]^2}{\si^2(\eta_s)}\; ds 
$$
for any fixed $0<t\le t_0$. By Cauchy-Schwarz combined with i), the mixed terms (\ref{term-2}) vanish in $P^\vth$-probability as $\nto$, uniformly over $[0,t_0]$. The lemma is proved. \halmos\\

For the next lemma and for the proof of theorem 1.2, we prepare some notation. Write $\calh$ for the family of (non-void) finite subsets of $\bbr$. Arrange elements of $H\in\calh$ in increasing order $h_1<h_2<\ldots<h_r$ with  $r=r(H)\in\bbn$, write $d=d(H)=\max_i |h_i|$, and let $\mathbb{A}^H$ denote the $r{\times}r$ matrix with $(i,j)$-entry
\beqq\label{for_my_limit-1}
(|h_i|\wedge|h_j|) \;\;\mbox{if ${\rm sgn}(h_i)={\rm sgn}(h_j)$} \quad,\quad 0 \;\;\mbox{if ${\rm sgn}(h_i)\neq{\rm sgn}(h_j)$} \;. 
\eeqq
This is the covariance matrix for $(\, \wt W (h_i) \,)_{1\le i\le r}$ in double sided Brownian motion $(\, \wt W (h) \,)_{h\in\bbr}$. 
For every $\vth\in\Theta$, we introduce an additional $r$-dimensional Brownian motion $\wh B$ which is independent of $B^\vth$. For $0<t_0<\infty$ fixed and $n\ge n_0(\vth,t_0,d)$, we consider continuous $r$-dimensional $P^\vth$-martingales
$$
M^{n,H,t_0}_t \;\;=\;\; \left(\, M^{n,h_i}_t \,\right)_{1\le i\le r} \quad,\quad t\ge 0\\[-2mm]
$$ 
defined by $dM^{n,H,t_0}_t = d\wh B_{t-t_0}$ for $t>t_0$ and
\beqq\label{for_my_limit-3}
M^{n,h_i}_t \;:=\;  \int_0^{tn} j^{h_i,n}_\vth(s)\, \frac{1}{\si(\eta_s)}\, dB^\vth_s \quad\mbox{for $t\le t_0$} \;, 
\eeqq
together with an $r$-dimensional Gaussian martingale (defined on some probability space)
$$ 
\wh W^{H,t_0} = \left( \wh W^{H,t_0}_t \right)_{t\ge 0}   \quad,\quad 
\wh W^{H,t_0}_t \;\;=:\;\; \left(\, \wh W^{h_i}_t \, \,\right)_{1\le i\le r} 
$$
defined by $d\wh W^{H,t_0}_t = d\wh B_{t-t_0}$ for $t>t_0$ and  
\beqq\label{for_my_limit-2}
\langle\, \wt W^{H,t_0} \,\rangle_t \;\;=\;\; \mathbb{A}^H\; t^2\; J(\vth,r,\rho)  \quad\mbox{for $t\le t_0$}
\eeqq
where $J(\vth,r,\rho)$ is the limiting constant in lemma 3.2. Processes $\wh W^{H,t_0}$ will be the limiting objects in the following lemma.\\

{\bf 3.4 Lemma: } Under {\bf (H0')}+{\bf (H1')}+{\bf (H2)}, for arbitrary $0<t_0<\infty$ fixed, we have weak convergence 
$$
M^{n,H,t_0} \quad\lra\quad \wh W^{H,t_0}    \quad\mbox{under $P^\vth$ as $\nto$}
$$
in the Skorohod path space $\bbd([0,\infty),\bbr^r)$. \\

{\bf Proof: } From 
the martingale convergence theorem ([JS 87, VIII.3.22]) combined with lemma 3.2. \halmos \\

Now we can prove theorem 1.2.\\

{\bf 3.5 Proof of theorem 1.2: } We assume {\bf (H0')}+{\bf (H1')}+{\bf (H2)}. We fix $\vth\in\Theta$ and $t_0>1$. For finite collections $H=\{h_1,\ldots,h_r\}\in\calh$, with $d=d(H)$, $r=r(H)$ and $n\ge n_0(\vth,t_0,d)$, we consider  
$$
\left(\, \log L_n^{ (\vth+h_i/n^2) / \vth } \,\right)_{1\le i\le r}
$$
which by proposition 3.3 is equivalent under $P^\vth$ as $\nto$ to 
$$
\left(\, \int_0^{n} j^{h_i,n}_\vth(s)\, \frac{1}{\si(\eta_s)}\, dB^\vth_s \;-\; \frac12 \int_0^{n} \left( j^{h_i,n}_\vth(s) \right)^2 \frac{1}{\si^2(\eta_s)}\, ds  \,\right)_{1\le i\le r} \;. 
$$
With notations of lemma 3.4 (since $t=1<t_0$) and  (\ref{for_my_limit-3}) this last object is 
$$
\left(\, M^{n,h_i}_1 \;-\; \frac12 \langle M^{n,h_i} \rangle_1 \,\right)_{1\le i\le r} \;.
$$ 
and converges by lemma 3.2 and lemma 3.4 weakly in $\bbr^r$  under $P^\vth$ as $\nto$ to 
$$
\left(\, \wh W^{h_i}_1 \;-\; \frac12 \langle \wh W^{h_i} \rangle_1 \,\right)_{1\le i\le r}  
$$
by definition of the process $\wh W^{H,t_0}$ above. By (\ref{for_my_limit-2}) and (\ref{for_my_limit-1}), this is equal in law to  
\beqq\label{mein_ziel}
\left(\, \wt W(\, h_i\, J(\vth,r,\rho)\,) \;-\; \frac12\,  |h_i|\, J(\vth,r,\rho) \,\right)_{1\le i\le r} 
\eeqq
for double-sided Brownian motion $(\wt W(h))_{h\in\bbr}$. 
\halmos\\

\newpage
{\Large\bf References} 



\vskip0.2cm
[BGT 87]\quad
Bingham, N., Goldie, C., Teugels, J.: Regular variation. Cambridge 1987. 

\vskip0.2cm
[B 05]\quad
Brandt, C: Partial reconstruction of the trajectories of a discretely observed branching diffusion with immigration and an application to inference. 
PhD thesis, Universit\"at Mainz 2005. 


\vskip0.2cm
[CLM 06]\quad
Castillo, I., L\'evy-Leduc, C., Matias, C.: 
Exact adaptive estimation of the shape of a periodic function with unknown period corrupted by white noise. 
Math.\ Methods Statist.\ {\bf 15}, 1--30 (2006).


\vskip0.2cm 
[D 10]\quad
Dachian, S.: 
On limiting likelihood ratio processes of som change-point type statistical models. 
J.~Statist.\ Plann.\ Inference {\bf 140}, 2682--2692 (2010). 

\vskip0.2cm 
[D 85]\quad
Davies, R.: Asymptotic inference when the amount of information is random. In: Le Cam, L., Olshen, R. (Eds): 
Proc. of the Berkeley Symposium  in honour of J. Neyman and J. Kiefer. Vol. II. Wadsworth 1985. 

\vskip0.2cm
[DFK 10]\quad
Dehling, H., Franke, B., Kott, T.: 
Drift estimation for a periodic mean reversion process. 
To appear in Statist.\ Inference Stoch.\ Proc.\ (2010). 

\vskip0.2cm 
[DP 84]\quad
Deshayes, J., Picard, D.: 
Lois asymptotiques des tests et estimateurs de rupture dans un mod\`ele statistique classique. 
Annales I.H.P.\ B {\bf 20}, 309--327 (1984). 

\vskip0.2cm
[G 79]\quad
Golubev, G.: 
Computation of efficiency of maximum-likelihood estimate when observing a discontinuous signal in white noise.  Problems Inform.\ Transmission  {\bf 15}, 61--69  (1979). 

\vskip0.2cm
[G 88]\quad
Golubev, G.: 
Estimating the period of a signal of unknown shape currupted by white noise. 
Problems Inform.\ Transmission {\bf 24}, 38--52 (1988). 


 

\vskip0.2cm 
[H 70]\quad 
H\'ajek, J.: A characterization theorem of limiting distributions for regular estimators.\\  
Zeitschr.\ Wahrscheinlichkeitstheor.\ Verw.\ Geb.\ {\bf 14}, 323--330, 1970. 
 
\vskip0.2cm 
[HS 67]\quad 
H\'ajek, J., Sid\'ak, Z:  Theory of rank tests. Academic Press 1967. 


\vskip0.2cm 
[HK 10]\quad 
H\"opfner, R., Kutoyants, Y.:  
Estimating discontinuous periodic signals in a time inhomogeneous diffusion. 
Statist.\ Inference Stoch.\ Proc.\ {\bf 13}, 193--230 (2010). 

\vskip0.2cm 
[HK 11]\quad 
H\"opfner, R., Kutoyants, Y.:  
On LAN for parametrized continuous periodic signals in a time inhomogeneous diffusion. 
To appear in Statistics$\&$Decisions (2011). 

\vskip0.2cm
[IH 81]\quad
Ibragimov, I., Has'minskii, R.: Statistical estimation. Springer 1981. 

\vskip0.2cm
[IW 89]\quad
Ikeda, N., Watanabe, S.: Stochastic differential equations and diffusion processes. 2nd ed.\ North-Holland/Kodansha 1989. 


\vskip0.2cm
[JS 87]\quad
Jacod, J., Shiryaev, A.: Limit theorems for stochastic processes. Springer 1987. 



\vskip0.2cm
[KS 91]\quad
Karatzas, J., Shreve, S.: Brownian motion and stochastic calculus. 2nd ed.\ Springer 1991. 

\vskip0.2cm
[KK 00]\quad
K\"uchler, U., Kutoyants, Y.: Delay estimation for some stationary diffusion-type processes. \\
Scand.\ J.\ Statist.\ {\bf 27}, 405--414 (2000). 

\vskip0.2cm
[K 98]\quad
Kutoyants, Y.: 
Statistical inference in spatial Poisson processes. 
Springer Lect.\ Notes Math.\ 134, Springer 1998. 

\vskip0.2cm
[K 04]\quad
Kutoyants, Y.: 
Statistical inference for ergodic diffusion processes. Springer 2004. 


\vskip0.2cm 
[L 68]\quad
Le Cam, L.: Th\'eorie asymptotique de la d\'ecision statistique. Presses de l'Universit\'e de Montr\'eal 1969. 
 
\vskip0.2cm 
[LY 90]\quad
Le Cam, L., Yang, G.: Asymptotics in statistics. Some basic concepts. 
Springer 1990.   (2nd Ed.\ Springer 2002).

\vskip0.2cm 
[LM 08]\quad
Liese, F., Miescke, K.: 
Statistical decision theory. Springer 2008. 

\vskip0.2cm
[LS 81]\quad
Liptser, R., Shiryaev, A.: Statistics of random processes. Vols.\ I+II,  Springer 1981, 2nd Ed.\ 2001. 



 
\vskip0.2cm
[P 94]\quad
Pfanzagl, J.: 
Parametric statistical inference. de Gruyter 1994. 




\vskip0.2cm
[RS 95]\quad
Rubin, H., Song, K.: Exact computation of the asymptotic efficiency of maximum likelihood estimators of a dicontinuous signal in a Gaussian white noise. Ann.\ Statist.\ {23}, 732--739 (1995). 


\vskip0.2cm
[S 85]\quad
Strasser, H.: Mathematical theory of  statistics. de Gruyter 1985. 


\vskip0.2cm
[V 98]\quad
van der Vaart, A.: 
Asymptotic Statistics. Cambridge 1998. 

\vskip0.2cm
[WM 95]\quad
Witting, H., M\"uller-Funk, U.: 
Mathematische Statistik II. Teubner 1995. 



\vfill
~\hfill{\bf 25.10.2010}\\ 

Yury A.\ Kutoyants\\
Laboratoire de Statistique et Processus, Universit\'e du Maine, F--72085 Le Mans Cedex 9\\ 
{\tt kutoyants@univ-lemans.fr}\\
{\tt http://www.univ-lemans.fr/sciences/statist/pages$\_$persos/kuto.html}

\vskip0.5cm
Reinhard H\"opfner\\
Institut f\"ur Mathematik, Universit\"at Mainz, D--55099 Mainz\\ 
{\tt hoepfner@mathematik.uni-mainz.de}\\
{\tt http://www.mathematik.uni-mainz.de/$\sim$hoepfner}

\end{document}